\documentclass[]{article}

\title{\bf More on a question of M. Newman on isomorphic subgroups of solvable groups}
\author{George Glauberman,\\
 Department of Mathematics, University of Chicago,\\
 5734 S. University Avenue,\\
 Chicago, IL 60637,\\ 
 \& \\
 Geoffrey R. Robinson,\\ Institute of Mathematics,\\ University of Aberdeen,\\ Aberdeen AB24 3UE}

\usepackage{latexsym,amssymb}
\begin{document}

\maketitle

\begin{abstract}
We make further remarks on a question of Moshe Newman, which asked whether it is the case that if $H$ and $K$ are isomorphic subgroups of a finite solvable group $G$ and $H$ is maximal in $G$, then $K$ is also maximal. This continues work begun in [1] by I.M. Isaacs and the second author.

\medskip
We prove here that if Newman's question has a negative answer for the triple $(G,H,K)$ (ie $H$ is maximal in $G$, but the isomorphic subgroup $K$ is not), then $p \leq 3$ and, for $q = 5-p,$ we have 
$$O_{q^{\prime}}(H) = O_{q^{\prime}}(K) = O_{q^{\prime}}(G)$$ and Newman's question also has a negative answer 
for the triple $(G^{\ast}, H^{\ast},K^{\ast}),$ where $G^{\ast} = G/ O_{q^{\prime}}(G)$, etc.. Furthermore, we prove that $G$ has a homomorphic image ${\bar G}$ such that Newman's question has a negative answer for the triple $({\bar G},{\bar H},{\bar K})$, while $F({\bar G}),F({\bar H})$ and $F({\bar K})$
are all $q$-groups, and $O_{\{2,3\}}({\bar H})$ involves ${\rm Qd}(q).$

\medskip
As an application, we prove that if $G$ is a finite solvable group such that $H$ and $K$ are isomorphic subgroups of $G$ with $H$ maximal and $K$ not maximal, with $[G:H] = [G:K]$ a power of the prime $p,$ then $p \leq 3$ and a Hall $\{2,3\}$-subgroup $L$ of $H$ necessarily involves $S_{3},A_{4}$ and a non-Abelian group of order $8$ (in fact, $L$ involves at least one of $S_{4}$ or ${\rm Qd}(3)).$ In particular, $L$ is neither $2$-closed nor $3$-closed.

\end{abstract}

\section{Introduction}

In this note, we make further remarks on a question of Moshe Newman, which asked whether it is the case that if $H$ and $K$ are isomorphic subgroups of a finite solvable group $G$ and $H$ is maximal in $G$, then $K$ is also maximal. This continues work begun in [1], where it was proved that this is true if $H$ either has Abelian Sylow $2$-subgroups or a Sylow tower, and it was noted (in Theorem 3 of that paper) that in a minimal counterexample $G$, we have $O_{p}(G) = 1$, where $[G:H]$ is a power of the prime $p$.

\section{Notation, Assumed Background, and Preliminary Results}

\medskip
Recall that a finite group $L$ is said to be involved in the finite group $G$ if there is a subgroup $H$ of $G$ and a normal subgroup $K$ of $G$ such that $H/K \cong L.$  We will make frequent use of the first lemma throughout.

\medskip
\noindent {\bf Lemma 1:} \emph{ Let $X$ be a finite solvable group, and $\pi$ be a set of primes. Let $T$ be a finite $\pi$-group and let $Y$ be a Hall $\pi$-subgroup of $X.$ Then $T$ is involved in $X$ if and only if $T$ is involved in $Y$.}

\medskip
\noindent {\bf Proof:} It is clear that $T$ is involved in $X$ if $T$ is involved in $Y$. We prove the opposite implication by induction on $|X|$. Suppose that $U,V$ are subgroups of $X$ with $V \lhd U$ and $U/V \cong T.$
Then $T$ is certainly involved in $U$. 

\medskip
If $|U| < |X|,$ then $T$ is involved in a Hall $\pi$-subgroup of $U$ by induction. But any Hall $\pi$-subgroup of $U$ is conjugate to a subgroup of $Y$, so that $T$ is involved in $Y$. Hence we may suppose that $U = X$ and that $V \lhd X$ with $X/V \cong T.$

\medskip
Let $M$ be a minimal normal subgroup of $X$ contained in $V.$ Then $M$ is either a $\pi$-group or a $\pi^{\prime}$-group. Now $T \cong (U/M)/(V/M)$ so that $T$ is involved in $X/M.$ By induction, $T$ is involved in the Hall $\pi$-subgroup $YM/M$ of $X/M.$ 

\medskip
If $M$ is a $\pi$-group then $Y/M$ is a Hall $\pi$-subgroup of $X/M$ and we are done since we already remarked that $T$ is involved in $Y/M.$ If $M$ is a $\pi^{\prime}$-group, then $YM/M \cong Y,$ so $T$ is involved in $Y,$ and the proof of the lemma is complete.

\medskip
Let $r$ be a prime. The group ${\rm Qd}(r)$ is the semi-direct product of the natural module for ${\rm SL}(2,r)$ 
with ${\rm SL}(2,r).$ Note that $|{\rm Qd}(r)| = r^{3}(r^{2}-1)$ and that ${\rm Qd}(2)$ is isomorphic to the symmetric group $S_{4}.$ Note also that ${\rm Qd}(3)$ is a $\{2,3\}$-group. We recall that ${\rm Qd}(r)$ is solvable if and only if $r < 5.$ In this note, the more elementary fact that $|{\rm Qd}(r)|$ has at least three prime divisors when $r > 3$ usually suffices for our purposes. For $r-1$ and $r+1$ can't both be powers of $2$ when $r$ is a prime greater than $3$.

\medskip
Using Lemma 1, we note that if ${\rm Qd}(r)$ is involved in a solvable group $X,$ then $r < 5$ and ${\rm Qd}(r)$ is involved in a Hall $\{2,3\}$-subgroup of $X.$ 

\medskip
Using Theorems of Stellmacher and Glauberman, and Lemma 1, for each finite $r$-group $R,$ there is a characteristic subgroup $W(R)$ of $R$ (which is non-trivial whenever $R$ is non-trivial) such that whenever $X$ is a finite solvable group with Sylow $r$-subgroup $R,$ then we have $X = O_{r^{\prime}}(X) N_{X}(W(R))$ and furthermore $O_{r^{\prime}}(X)W(R)$ is characteristic in $X,$ unless, perhaps, $r <5$ and ${\rm Qd}(r)$ is involved in a Hall $\{2,3\}$-subgroup of $X$. When $r$ is odd, we may take $W(R) = ZJ(R).$

\medskip
When $n$ is an integer, we let $\pi(n)$ denote the set of prime divisors of $n$.
We will make frequent use of the following easy lemma:

\medskip
\noindent {\bf Lemma 2:} \emph{ Let $X$ be a finite solvable group and $\pi$ be a set of primes.
Let $Y$ be a subgroup of $X$ with $\pi([X:Y]) \subseteq \pi$. Then $O_{\pi}(Y) \leq O_{\pi}(X).$ }

\medskip
\noindent {\bf Proof:} Let $U$ be a Hall $\pi$-subgroup of $X$ containing $O_{\pi}(Y).$  Since $Y$ contains a Hall $\pi^{\prime}$-subgroup of $X,$ we have $X = YU = UY.$ Now we have 
$$O_{\pi}(X) = \cap_{x \in X} U^{x} = \cap_{y \in Y} U^{y} \geq O_{\pi}(Y),$$ as claimed.

\medskip
\noindent {\bf Corollary 3:} \emph{Let $X$ be a finite solvable group and $Y$ be a subgroup of $X$ whose index is a power of the prime $s.$ Let $r$ be a prime divisor of $X$ different from $s$. Then $O_{\{r,s\}}(Y)$ and $O_{\{r,s\}}(X)$ have a common Sylow $r$-subgroup.}

\medskip
\noindent {\bf Proof:} By Lemma 2, applied with $\pi = \{r,s\}$, we have $O_{\{r,s\}}(Y) \leq O_{\{r,s\}}(X).$ Note that $Y$ contains a Sylow $r$-subgroup of $X,$ so that $Y$ certainly contains a Sylow $r$-subgroup $R$ of $O_{\{r,s\}}(X).$  Now $$R \cap Y  \leq O_{\{r,s\}}(X) \cap Y \leq O_{\{r,s\}}(Y),$$  so that 
$O_{\{r,s\}}(Y)$ contains a Sylow $r$-subgroup of $O_{\{r,s\}}(X).$ On the other hand, since
$O_{\{r,s\}}(Y) \leq O_{\{r,s\}}(X),$ a Sylow $r$-subgroup of $O_{\{r,s\}}(Y)$ is contained 
in some Sylow $r$-subgroup of $O_{\{r,s\}}(X).$

\section{ Statement and Proof of Theorem A}

\medskip
\noindent {\bf Theorem A:} \emph{ Let $H$ be a maximal subgroup of the finite solvable group $G$  and suppose that $[G:H] = p^{a}$ where $p$ is a prime and $a$ is a positive integer. Let $K$ be a subgroup of $G$ which is isomorphic to $H.$ Suppose that $K$ is not maximal in $G.$}

\medskip
\emph{Then $p \leq 3,$ and, for $q = 5-p,$ we have $$O_{q^{\prime}}(H) = O_{q^{\prime}}(G) = O_{q^{\prime}}(K)$$
and, for $G^{\ast} = G/O_{q^{\prime}}(G),$ etc., $H^{\ast}$ and $K^{\ast}$ are isomorphic subgroups 
of $G^{\ast}$ with $H^{\ast}$ maximal and $K^{\ast}$ not maximal.}

\medskip
\noindent {\bf Proof:} Let $\phi: H \to K$ be an isomorphism.

\medskip
We proceed by induction on $|G|$. Note that if $N$ is a $\phi$-invariant normal subgroup of $H$ then 
$$N = N\phi \lhd H\phi = K$$ and then $N \lhd \langle H,K \rangle = G.$ Then $\phi$ induces an isomorphism 
between $H/N$ and $K/N$, and $H/N$ is maximal in $G/N$, but $K/N$ is not maximal in $G/N$.

\medskip
Using Theorem 3 of [1] and Lemma 2, we may conclude that $O_{p}(G) = 1.$ For if not, then we have 
$O_{p}(G) \subseteq O_{p}(G) \cap H$ by the former result, and both $O_{p}(H) \leq O_{p}(G)$ and  
$O_{p}(K) \leq O_{p}(G)$ by the latter. Then we have $$O_{p}(G) = O_{p}(H)$$ and
$$O_{p}(H)\phi = O_{p}(K) \leq O_{p}(G).$$ Since $$|O_{p}(G)| = |O_{p}(H)| = |O_{p}(K)|,$$
we have $$O_{p}(G)\phi = O_{p}(G),$$ so we may apply the argument above with $O_{p}(G)$ in the role of $N$.

\medskip
We may suppose by induction that the theorem holds for the triple $(G/N,H/N,K/N).$ Then $p \leq 3$
since $[G/N:H/N] = [G:H],$ and we note that $$G/O_{q^{\prime}}(G) \cong (G/N)/O_{q^{\prime}}(G/N)$$ for $q = 5-p$. Hence the theorem holds for $G$ in this case. Thus we may suppose that $O_{p}(G) = 1.$

\medskip
Since $O_{p}(G) = 1$ and (by Lemma 2) $O_{p}(H) \leq O_{p}(G)$, we see that $F(H)$ is a $p^{\prime}$-group, as is the isomorphic group $F(K)$. Also, there is a prime $r \neq p$ such that $O_{r}(H) \neq 1.$ 

\medskip
Suppose first that $r > 3.$ Then by Lemma 1, we have $$O_{\{r,p\}}(H) \leq O_{\{r,p\}}(G)$$ and likewise,
$$O_{\{r,p\}}(K) \leq O_{\{r,p\}}(G).$$ Furthermore, by Corollary 3 we may suppose (possibly after replacing $K$ by a conjugate) that $O_{\{r,p\}}(H),O_{\{r,p\}}(G)$ and $O_{\{r,p\}}(K)$ all have a common Sylow $r$-subgroup, say $S$. 

\medskip
Since $O_{p}(H) = O_{p}(K) = O_{p}(G) = 1,$ we see that $W(S)$ is characteristic in each of 
$O_{\{r,p\}}(H),O_{\{r,p\}}(K)$ and $O_{\{r,p\}}(G),$ where $W(S)$ is the Glauberman-Stellmacher characteristic subgroup of $S$ (for note that $O_{\{r,p\}}(G)$ does not involve ${\rm Qd}(2)$ or ${\rm Qd}(3)$ since $r > 3,$
and likewise for $O_{\{r,p\}}(H)$ and $O_{\{r,p\}}(K))$. Thus $W(S) \lhd G.$ Furthermore, $$O_{\{r,p\}}(H)\phi = O_{\{r,p\}}(K)$$ and composing $\phi$ with an inner automorphism of $K$ if necessary, we may suppose that $S\phi = S,$ in which case $W(S)\phi = W(S).$ Then $\phi$ induces an isomorphism between $H/W(S)$ and $K/W(S).$  

\medskip
By induction, the theorem holds for $G/W(S).$ In particular, $p \leq 3,$ and, setting $q = 5-p,$ we see that $W(S)$ is a $q^{\prime}$-group, and that the Theorem therefore holds for $G.$

\medskip
Hence we may suppose that $F(H)$ is a $q$-group for some prime $q \leq 3.$\\ 
If $p \neq 5-q,$ we may argue as above that $W(T)$ is normal in each of\\ $O_{\{q,p\}}(H),O_{\{q,p\}}(G)$ and $O_{\{q,p\}}(K),$ where $T$ is a common Sylow $q$-subgroup of $O_{\{q,p\}}(H),O_{\{q,p\}}(G)$ and $O_{\{q,p\}}(K).$ This time, $O_{\{q,p\}}(G)$ does not involve ${\rm Qd}(2)$ or ${\rm Qd}(3)$ since $p > 3.$
But that leads to a contradiction, since the theorem holds for $G/W(T)$ and then we see that 
$[G:H] = [G/W(T):H/W(T)]$ is either a power of $2$ or power of $3,$ so $p \leq 3,$ contrary to current assumptions.

\medskip
Now we may suppose that $p \leq 3$ and that $F(H)$ is a $q$-group, where $q = 5-p.$ Since $[G:H]$ is a power of $p$ and $O_{p}(G) = 1,$ we have $F(G) \leq H,$ so that $F(G) \leq F(H)$ and $F(G)$ is a $q$-group. 
Since $H \cong K,$ we also see that $F(K)$ is a $q$-group. The proof of Theorem A is complete.

\section{Statement and Proof of Theorem B and some consequences}

\medskip
\medskip
\noindent {\bf Theorem B:} \emph{ Let $H$ be a maximal subgroup of the finite solvable group $G$ and suppose that $[G:H] = p^{a}$ where $p \leq 3$ is a prime and $a$ is a positive integer. Let $K$ be a subgroup of $G$ which is isomorphic to $H.$ Suppose that $K$ is not maximal in $G$ and that $F(H),F(K)$ and $F(G)$ are all
$q$-groups, where $q = 5-p.$  Let $Q$ be a Sylow $q$-subgroup of $H$.}  

\medskip
\emph{ Then $G$ has a homomorphic image $G^{\ast}$ such that $H^{\ast}$ and $K^{\ast}$ (the respective images of $H$ and $K$) are isomorphic subgroups of $G^{\ast}$ with $H^{\ast}$ maximal and $K^{\ast}$ not maximal, and with $F(H^{\ast}),F(K^{\ast})$ and $F(G^{\ast})$ all $q$-groups. Furthermore, $O_{\{2,3\}}(H^{\ast})$ involves ${\rm Qd}(q)$ and no non-identity characteristic subgroup of $Q^{\ast}$ is normal in $H^{\ast}.$ } 

\medskip
\noindent {\bf Proof:} Let $\phi:H \to K$ be the isomorphism as before. We may, and do, suppose that $Q\phi = Q,$
composing $\phi$ with an inner automorphism of $K$ if necessary. 

\medskip
Using Theorem A, we may suppose that there is no non-trivial $\phi$-invariant normal subgroup of $H.$
For otherwise, if $1 \neq N \lhd H$ is $\phi$-invariant and chosen of maximal order subject to these conditions, 
then $N = N \phi \lhd H\phi = K,$ so $N \lhd \langle H,K \rangle = G$ and $\phi$ induces an isomorphism between $H/N$ and $K/N$. By Theorem A, we have $$O_{q^{\prime}}(H/N) = O_{q^{\prime}}(G/N) = O_{q^{\prime}}(K/N)= 1$$ by the maximal choice of $N$ (hence $F(H/N),F(K/N)$ and $F(G/N)$ are all $q$-groups) and the theorem holds for $G/N$ by induction, so it holds for $G.$

\medskip
Suppose now that $O_{\{2,3\}}(H)$ does not involve ${\rm Qd}(q).$ Then as before, $Q$ contains a common Sylow $q$-subgroup $T$ of $O_{\{2,3\}}(H), O_{\{2,3\}}(K)$ and $O_{\{2,3\}}(G)$ such that $W(T) \lhd H$ 
and $W(T)\phi = W(T),$ contradicting the fact that there is no non-trivial $\phi$-invariant normal subgroup of $H$. Hence $O_{\{2,3\}}(H)$ involves ${\rm Qd}(q).$

\medskip
If there is a non-identity characteristic subgroup $S$ of $Q$ which is normal in $H,$ then we have $S = S\phi$ since $Q = Q\phi,$ again contradicting the fact that  no non-trivial $\phi$-invariant normal subgroup of $H$. The proof of Theorem B is complete.

\medskip
We may combine Theorems A and B to deduce:

\medskip
\noindent {\bf Corollary C:} \emph{  Let $G$ be a finite solvable group and 
$H,K$ be isomorphic subgroups of $G$ such that $H$ is maximal but $K$ is not, and with $[G:H] = p^{a}$ for some prime $p$ and positive integer $a$. Suppose (with no loss of generality) that $H \cap K$ contains a Hall $p^{\prime}$-subgroup of $G$ (this is just a matter of replacing $K$ by a conjugate if necessary). Then $p \leq 3$. Let $\phi: H \to K$ be an isomorphism, chosen so that $Q\phi = Q$ for a Sylow $q$-subgroup $Q$ of $H,$ where $q = 5-p$ (this can be achieved by composing $\phi$ with an inner automorphism of $K$ if necessary). Let $C$ be the unique maximal $\phi$-invariant normal subgroup of $G$ which is contained in $H \cap K,$  and let $G^{\ast} = G/C,$ etc.. Then $F(H^{\ast}),F(K^{\ast})$ and $F(G^{\ast})$  are all $q$-groups. Also, ${\rm Qd}(q)$ is involved in $O_{\{2,3\}}(H^{\ast}),$ and no non-identity characteristic subgroup of $Q^{\ast}$ is normal in $H^{\ast}.$ }

\medskip
\noindent {\bf Proof:} The proofs of Theorems A and B show that unless $G$ itself can play the role of $G^{\ast},$ there is a non-identity $\phi$-invariant normal subgroup $N$ of $G$ contained in both $H$ and $K$ such that  $H/N$ and $K/N$ are isomorphic with $H/N$ maximal in $G/N$ and $K/N$ not maximal. In that case, since
$N \leq C,$ we may suppose that the Theorem holds in $G/N$ by induction, and then it holds in $G.$ 

\medskip
\noindent {\bf Corollary D:} \emph{Let $G$ be a finite solvable group containing isomorphic subgroups $H$ and $K$ such that $H$ is maximal in $G$ but $K$ is not. Then $[G:H] = [G:K]$ is a power of a prime $p \leq 3$ and a Hall $\{2,3\}$-subgroup $L$ of $H$ involves ${\rm Qd}(q),$ where $q = 5-p.$ In particular, $L$ involves $S_{3},A_{4}$ and a non-Abelian group of order $8,$ so that $L$ is neither $2$-closed nor $2$-nilpotent.} 

\medskip
\begin{center}
{\bf References}
\end{center}

\medskip
\noindent [1] Isaacs, I.M. \& Robinson, G.R., \emph{Isomorphic subgroups of solvable groups}, Proc Amer. Math. Soc
{\bf 143} (2015) 3371-3376.

\end{document}